\documentclass[12pt]{amsproc}
\usepackage[cp1251]{inputenc}
\usepackage[T2A]{fontenc} 
\usepackage[russian]{babel} 
\usepackage{amsbsy,amsmath,amssymb,amsthm,graphicx,color} 
\usepackage{amscd}
\def\le{\leqslant}
\def\ge{\geqslant}

\newtheorem{prop}{Предложение}

\theoremstyle{definition}

\theoremstyle{remark}

\begin {document}
\unitlength=1mm
\title[Производящие функции Харера-Цагира]
{Производящие функции Харера-Цагира, цикловой индекс Редфилда-Пойа и полулинейные сравнения Коэна}
\author{Г. Г. Ильюта}
\email{ilgena@rambler.ru}
\address{}
\thanks{Работа поддержана грантом РФФИ-16-01-00409}


\begin{abstract}
Производящие функции Харера-Цагира для эйлеровых характеристик пространств модулей кривых содержат $n$-ожерельные многочлены. Разложения Тейлора для этих многочленов зависят от количеств решений полулинейных сравнений Коэна.

Harer-Zagier generating functions for Euler characteristics of moduli spaces of curves contain $n$-necklace polynomials. Taylor expansions for these polynomials depend on numbers of solutions of Cohen semilinear congruences.
\end{abstract}

\maketitle

  В производящих функциях для эйлеровых характеристик пространств модулей кривых \cite{6}, Th. 4', Th. 5'; \cite{3}, Th. 4.5, Remark 4.6, появляются многочлены $\beta_{k,d}(t)$. Мы докажем в п. 1, что они отличаются простым преобразованием от обобщённых цикловых индексов Редфилда-Пойа, отвечающих регулярному представлению циклической группы $C_k$ порядка $k$ и её неприводимым характерам. Обобщённые цикловые индексы известны также как симметрические функции Шура конечной группы перестановок (характеры общей линейной группы) \cite{10} или как образы при отображении Фробениуса индуцированных характеров симметрической группы \cite{11}, p. 395. Многочлены $\beta_{k,d}(t)$ были введены для алгебраического упрощения производящих функций для эйлеровых характеристик. Начальные варианты этих производящих функций содержали количества решений некоторых линейных сравнений (они появляются при гомоморфизме в циклическую группу как образы соотношений в исходной группе), которые затем были "упакованы"\: в многочлены $(\beta_{k,d}(t))^s$, $s\in\mathbb Z_{>0}$. Одна из частей равенства, приводящего к упрощению производящих функций \cite{6}, р. 482, определяется количествами решений линейных сравнений -- в п. 3 мы свяжем другую часть с количествами решений полулинейных сравнений Коэна. Возможно, поиск прямой связи между полулинейными сравнениями Коэна и комбинаторной топологией из \cite{6} и \cite{3} мог бы привести к сжатию информации не только на алгебраическом уровне производящих функций, но и на уровне комбинаторной топологии, что позволило бы упростить доказательства в \cite{6} и \cite{3}. Полулинейные сравнения Коэна можно рассматривать как семейства линейных сравнений в следующем смысле. Если $(x_1,\dots,x_m;y_1,\dots,y_m)$ -- решение сравнения Коэна (7) для $r=1$, то $(x_1,\dots,x_m)$ -- решение линейного сравнения, коэффициенты которого зависят от $(y_1,\dots,y_m)$ (формула для количества решений линейного сравнения с произвольными коэффициентами имеется в \cite{2}, p. 138). Такая интерпретация может быть полезной для комбинаторных (биективных) доказательств формул из п. 3.

  Для $t\in\mathbb Z_{>0}$ специализация обобщённого циклового индекса в правой части формулы (4) совпадает с размерностью отвечающего неприводимому характеру класса симметрии тензоров, причём, это верно для обобщённого циклового индекса любого неприводимого характера любой конечной группы перестановок \cite{10}, p. 226. Поэтому для регулярного представления циклической группы формулы из п. 3 связывают размерности классов симметрии тензоров, отвечающие неприводимым характерам циклической группы, с количествами решений полулинейных сравнений Коэна. Интересно было бы найти категорификацию этих формул.

  В п. 1 в качестве промежуточных объектов между обобщёнными цикловыми индексами Редфилда-Пойа и многочленами Харера-Цагира $\beta_{k,d}(t)$ вводятся $n$-ожерельные многочлены $M(t;n,k)$ (они определяются классическими суммами Рамануджана $c(n,k)$). Выделение этих многочленов мотивируется связями ожерельных многочленов $M(t;1,k)$ с различными разделами математики (обзор имеется в \cite{7}) -- возникает вопрос о связях этих разделов с пространствами модулей кривых. Другие обобщения многочленов  $M(t;1,k)$ \cite{7} определяют эйлеровы характеристики пространств модулей неприводимых многочленов от нескольких переменных над $\mathbb R$ и над $\mathbb C$ -- вопрос о связи этих пространств с пространствами модулей кривых также остаётся открытым. Для $t\in\mathbb Z_{>0}$ производящий многочлен для чисел $M(t;n,k)$
$$
\sum_{n=1}^kM(t;n,k)q^n
$$ 
является специализацией симметрической функции Краскиевича-Веймана \cite{9}, p. 9, \cite{1}, p. 9 (полагаем $t$ переменных этой функции равными $1$, а остальные -- равными $0$).

  В п. 2 мы покажем как зависит от многочленов $M(t;n,k)$ производящая функция для эйлеровых характеристик $e(\Gamma_g^1)$ групп классов отображений $\Gamma_g^1$ кривых рода $g$ с одной отмеченной точкой \cite{6}, Th. 4'. Используется равенство (3), связывающее многочлены $M(t;n,k)$ и $\beta_{k,d}(t)$. Также из этого равенства вытекает, что появляющиеся в \cite{5}, p. 447, многочлены являются обобщениями для любой конечной группы многочленов Харера-Цагира $\beta_{k,d}(t)$ (с заменой $t^a\to t^{k-a}$ для всех $a$). Роль сумм Рамануджана играют суммы значений неприводимого характера группы на всех элементах этой группы, имеющих фиксированный порядок.

  В п. 3 доказаны формулы, связывающие количества решений $Q_r(n,k,m)$ полулинейных сравнений Коэна с многочленами $M_r(t;n,k)$ (они определяются $r$-суммами Рамануджана $c_r(n,k)$ из \cite{4}). С многочленами Харера-Цагира $\beta_{k,d}(t)$ связан частный случай $Q_1(n,k,s)$ и $M_1(t;n,k)=M(t;n,k)$, но многочлены $M_r(t;n,k)$ позволяют использовать числа $Q_r(n,k,m)$ для всех $r$. Доказанные в п. 3 формулы вытекают из равенства
$$
\left.\left(t\frac{d}{dt}\right)^mM_r(t;n,k)\right|_{t=1}=\frac{Q_r(n,k,m)}{k^{rm}}.                \eqno (1)
$$
Соотношения между числами $Q_r(n,k,m)$ и многочленами $M_r(t;n,k)$ можно разными способами представить как соотношения между производящими функциями для них, мы рассмотрим один такой пример. В \cite{4}, p. 548, получена формула для производящего ряда Дирихле чисел $Q_r(n,k,m)$, мы докажем аналогичную формулу для многочленов $M_r(t;n,k)$. Для соотношений между этими рядами Дирихле формула (1) сводится к равенству
$$
\left.\left(t\frac{d}{dt}\right)^mLi_s(t)\right|_{t=1}=\zeta(s-m),             
$$
где $\zeta(s)=Li_s(1)$ -- дзета-функция Римана и $Li_s(t)$ -- полилогарифм,
$$
Li_s(t):=\sum_{k\ge 1}\frac{t^k}{k^s}.
$$

  1. Многочлены Харера-Цагира и $n$-ожерельные многочлены. Для $n,k\in\mathbb Z_{>0}$ определим $n$-ожерельные многочлены $M(t;n,k)$ формулой
$$
M(t;n,k):=\frac{1}{k}\sum_{d|k}c(n,k/d)t^d,
$$
где $c(n,k)$ -- суммы Рамануджана,
$$
c(n,k):=\sum_{\substack{1\le m\le k\\ (m,k)=1}}\epsilon_k^{mn},
$$
$\epsilon_k:=e^\frac{2\pi i}{k}$, $(n,k)$ -- наибольший общий делитель чисел $n$ и $k$. В частности, $c(0,k)=\phi(k)$ -- функция Эйлера, $c(1,k)=\mu(k)$ -- функция Мёбиуса.

  Многочлены $\beta_{k,d}(t)$ определяются формулой \cite{6}, p. 482,
$$
\beta_{k,d}(t):=\sum_{r=1}^{k-1}\epsilon_d^rt^{k-(k,r)}=\sum_{l|k,l< k}c(k,l,d)t^{k-l}, \qquad d|k,
$$
где для $l|k$ и $d|k$
$$
c(k,l,d):=\sum_{\substack{1\le m\le k\\ (m,k)=l}}\epsilon_d^m.
$$

  Для характера $\chi$ группы перестановок $H$ на множестве из $v$ элементов обобщённый цикловой индекс $Z_H^\chi(t_1,\dots,t_v)$ определяется формулой
$$
Z_H^\chi(t_1,\dots,t_v):=\frac{1}{|H|}\sum_{h\in H}\chi(h)\prod_{j=1}^vt_j^{c_j(h)},
$$
где $c_j(h)$ -- количество циклов длины $j$ в перестановке $h$. Для конечной группы $G$ обозначим через $R(G)$ её образ при регулярном представлении. Тогда $R(G)$ -- группа перестановок на множестве $G$. Через $\chi_n$ обозначим неприводимый характер циклической группы $C_k$, значение которого на образующей группы равно $\epsilon_k^n$.

\begin{prop}\label{prop1} 
$$
c(k,l,d)=c(k/d,k/l),                  \eqno (2)
$$
$$
\beta_{k,d}(t)=kt^kM(1/t;k/d,k)-1,    \eqno (3)              
$$
$$
M(t;n,k)=Z_{R(C_k)}^{\chi_n}(t,\dots,t).  \eqno (4)
$$
\end{prop}

  Доказательство. Полагая $a:=m/l$, получим
$$
c(k,l,d)=\sum_{\substack{1\le m\le k\\ (m,k)=l}}\epsilon_d^m=\sum_{\substack{1\le m\le k\\ (m/l,k/l)=1}}\epsilon_{k/l}^{mk/(dl)}.
$$
$$
=\sum_{\substack{1\le a\le k/l\\ (a,k/l)=1}}\epsilon_{k/l}^{ak/d}=c(k/d,k/l).
$$
Из этого равенства следует формула (3). При регулярном представлении группы $G$ элемент $g\in G$ (обозначим через $o(g)$ его порядок) действует как циклическая перестановка на каждом классе смежности порождённой им подгруппы, а значит является произведением $|G|/o(g)$ циклов длины $o(g)$. Поэтому
$$
Z_{R(G)}^{\chi}(t_1,\dots,t_{|G|})=\frac{1}{|G|}\sum_{d||G|}t_{|G|/d}^d\sum_{o(g)=|G|/d}\chi(g).
$$
В частности, для $G=C_k$ из определения сумм Рамануджана следует равенство
$$
Z_{R(C_k)}^{\chi_n}(t_1,\dots,t_k)=\frac{1}{k}\sum_{d|k}c(n,k/d)t_{k/d}^d. \blacksquare
$$

2. Производящие функции Харера-Цагира. Согласно \cite{6}, Th. 4', производящая функция для эйлеровых характеристик $e(\Gamma_g^1)$ групп классов отображений $\Gamma_g^1$ кривых рода $g$ с одной отмеченной точкой имеет вид
$$
\sum_{g\ge 1}e(\Gamma_g^1)t^{2g-1}=\sum_{k\ge 1}\frac{\phi(k)}{k}\sum_{d|k}\mu(d)\Phi^1(\beta_{k,d}(t),kt^k),    \eqno (5)
$$
где
$$
\Phi^1(X,Y):=\sum_{s\ge 2}\frac{(-1)^sX^s}{s(s-1)Y}
-\sum_{\substack{h\ge 1\\s\ge 0}}\binom{s+2h-2}{s}\frac{(-1)^sB_{2h}X^sY^{2h-1}}{2h}
$$
$$
=\frac{1}{Y}((1+X)\log(1+X)-X)+\mathfrak B\left(\frac{Y}{1+X}\right),
$$
$$
\mathfrak B(T)=-\sum_{h\ge 1}\frac{B_{2h}}{2h}T^{2h-1}\in\mathbb Q[[T]],
$$
$B_{2h}$ -- числа Бернулли. С помощью Предложения 1 формула (5) преобразуется к следующему виду.

\begin{prop}\label{prop2} 
$$
\sum_{g\ge 1}e(\Gamma_g^1)t^{2g-1}=\sum_{k\ge 1}\frac{\phi(k)}{k}\sum_{d|k}\mu(d)\bar\Phi^1(M(1/t;k/d,k),kt^k), 
$$
где
$$
\bar\Phi^1(Z,Y)=\Phi^1(YZ-1,Y)=Z\log(YZ)-Z+1/Y+\mathfrak B(1/Z). \blacksquare
$$ 
\end{prop}
  
  Аналогично зависит от многочленов $M(t;n,k)$ производящая функция для эйлеровых характеристик $e(\Gamma_g)$ групп классов отображений $\Gamma_g$ кривых рода $g$ без отмеченных точек \cite{6}, Th. 5',
$$
\sum_{g\ge 1}e(\Gamma_g)t^{2g-2}=\sum_{k\ge 1}\sum_{m,d|k}\frac{\phi(d)\mu(m)}{m^2}\Phi\left(\beta_{k/m,d/(d,m)}(t^m),\frac{kt^k}{m}\right),
$$
где
$$
\Phi(X,Y):=\sum_{s\ge 3}\frac{(-1)^{s-1}X^s}{s(s-1)(s-2)}+\sum_{s\ge 1}\frac{(-1)^sX^sY^2}{12s}
$$
$$
+\sum_{h\ge 2}\frac{B_{2h}}{2h(2h-2)}\left(\frac{Y}{1+X}\right)^{2h-2}.
$$

  3. Полулинейные сравнения Коэна. Для делителей $l_1,\dots,l_s$ числа $k$ пусть $N_k(b;l_1,\dots,l_s)$ -- количествo решений сравнения
$$
x_1+\dots+x_s=b\mod k,\quad (x_i,k)=l_i,\quad i=1,\dots,s.
$$ 
Согласно \cite{2}, p. 137,
$$
N_k(b;l_1,\dots,l_s)=\frac{1}{k}\sum_{d|k}c(b,d)\prod_{i=1}^sc(k/d,k/l_i).
$$
В \cite{6}, p. 479, появляются частные случаи чисел $N_k(b;l_1,\dots,l_s)$ и они удаляются из производящих функций \cite{6}, p. 481, с помощью соответствующих частных случаев равенства (мы используем формулы (2) и (3))
$$
\sum_{\substack{l_1,\dots,l_s|k\\l_i\ne k,i=1,\dots,s}}N_k(b;l_1,\dots,l_s)t^{ks-\sum l_i}
$$
$$
=\frac{1}{k}\sum_{d|k}c(b,d)\sum_{\substack{l_1,\dots,l_s|k\\l_i\ne k,i=1,\dots,s}}c(k/d,k/l_1)t^{k-l_1}\dots c(k/d,k/l_s)t^{k-l_s} 
$$
$$
=\frac{1}{k}\sum_{d|k}c(b,d)(\beta_{k,d}(t))^s.     
$$ 
Аналогично для многочленов $M(t;n,k)$ имеем
$$
\sum_{l_1,\dots,l_s|k}N_k(b;l_1,\dots,l_s)t^{\sum l_i}
$$
$$
=\frac{1}{k}\sum_{d|k}c(b,d)\sum_{l_1,\dots,l_s|k}c(k/d,k/l_1)t^{l_1}\dots c(k/d,k/l_s)t^{l_s}
$$
$$
=\frac{1}{k}\sum_{d|k}c(b,d)(kM(t;k/d,k))^s.      \eqno (6)
$$ 
  
  Пусть $Q_r(n,k,m)$ -- количество решений
$$
x_i\mod k,\quad y_i\mod k^r,\quad i=1,\dots,m,
$$
полулинейного сравнения Коэна 
$$
a_1x_1^ry_1+\dots+a_sx_m^ry_m=n\mod k^r,                 \eqno (7)
$$
где $(a_i,k)=1$ для всех $i$. Согласно \cite{4}, p. 547,
$$
\frac{Q_r(n,k,m)}{k^{rm}}=\frac{1}{k^r}\sum_{d|k}c_r(n,k/d)d^m.  \eqno (8)
$$ 

  В Предложении 3 мы свяжем числа $Q_r(n,k,m)$ с более общими многочленами
$$
M_r(t;n,k):=\frac{1}{k^r}\sum_{d|k}c_r(n,k/d)t^d,
$$
где $c_r(n,k)$ -- $r$-суммы Рамануджана \cite{4},
$$
c_r(n,k):=\sum_{\substack{1\le m\le k^r\\(m,k^r)_r=1}}\epsilon_{k^r}^{mn},
$$
$(a,b)_r$ -- наибольший общий делитель чисел $a$ и $b$, являющийся $r$-й степенью. В частности, $c_1(n,k)=c(n,k)$ и $M_1(t;n,k)=M(t;n,k)$. 

   Числа Стирлинга первого рода $s(n,m)$ определяются равенством
$$
\prod_{i=0}^{n-1}(t-i)=\sum_{m=0}^ns(n,m)t^m.
$$

  Пусть
$$
\delta_{k^r|n}:=M_r(1;n,k)=\frac{1}{k^r}\sum_{d|k}c_r(n,k/d).
$$
Согласно \cite{4}, p. 546, сумма в правой части равна $1$, если $k^r|n$, и равна $0$ в других случаях.

\begin{prop}\label{prop3} 
$$
\left.\frac{d^l}{dt^l}M_r(t;n,k)\right|_{t=1}=\sum_{m=1}^ls(l,m)\frac{Q_r(n,k,m)}{k^{rm}},\quad l\ge 1,           \eqno (9)
$$
$$
M_r(t;n,k)=\delta_{k^r|n}+\sum_{l=1}^k\frac{(t-1)^l}{l!}\sum_{m=1}^ls(l,m)\frac{Q_r(n,k,m)}{k^{rm}},              \eqno (10)
$$
$$
M_r(e^\lambda;n,k)=\delta_{k^r|n}+\sum_{l=1}^{\infty}\frac{\lambda^l}{l!}\frac{Q_r(n,k,l)}{k^{rl}},               \eqno (11)
$$
в частности, для $r=1$
$$
\left.\frac{d^l}{dt^l}M(t;n,k)\right|_{t=1}=\sum_{m=1}^ls(l,m)\frac{Q_1(n,k,m)}{k^{m}}, \quad l\ge 1,          
$$
$$
M(t;n,k)=\delta_{k|n}+\sum_{l=1}^k\frac{(t-1)^l}{l!}\sum_{m=1}^ls(l,m)\frac{Q_1(n,k,m)}{k^m},
$$
$$
M(e^\lambda;n,k)=\delta_{k|n}+\sum_{l=0}^{\infty}\frac{\lambda^l}{l!}\frac{Q_1(n,k,l)}{k^l}.
$$
\end{prop}

  Доказательство.
$$
\left.\left(t\frac{d}{dt}\right)^mM_r(t;n,k)\right|_{t=1}=\frac{1}{k^r}\left.\left(t\frac{d}{dt}\right)^m\sum_{d|k}c_r(n,k/d)t^d\right|_{t=1}
$$
$$
=\left.\frac{1}{k^r}\sum_{d|k}c_r(n,k/d)d^mt^d\right|_{t=1}=\frac{Q_r(n,k,m)}{k^{rm}}.
$$
Используя равенство \cite{8}, p. 4,
$$
t^l\left(\frac{d}{dt}\right)^l=\sum_{m=1}^ls(l,m)\left(t\frac{d}{dt}\right)^m,
$$
получим формулу (9). Формулы (10) и (11) представляют собой разложения Тейлора с учётом равенства 
$$
\left.\left(\frac{d}{d\lambda}\right)^lM_r(e^\lambda;n,k)\right|_{\lambda=0}=\left.\left(t\frac{d}{dt}\right)^lM_r(t;n,k)\right|_{t=1}=\frac{Q_r(n,k,l)}{k^{rl}}.\blacksquare
$$

  Отметим ещё один способ связать числа $N_k(b;l_1,\dots,l_s)$ и $Q_1(n,k,s)$. Заменяя в правой части формулы (6) $t^d$ на $d^t$ для всех $d$, получим следующий её аналог 
$$
\sum_{l_1,\dots,l_s|k}N_k(b;l_1,\dots,l_s)(\prod l_i)^t
=\frac{1}{k}\sum_{d|k}c(b,d)(\sum_{\delta|k}c(k/d,k/\delta)\delta^t)^s.
$$
Значения левой части этой формулы в точках $t\in\mathbb Z_{\ge 0}$, а значит и коэффициенты интерполяционного ряда Ньютона этой функции, определяются числами $Q_1(n,k,s)$ (это следует из формулы (8)). Напомним, что ряд Ньютона для функции $h(t)$ определяется по любой последовательности $a_0,a_1,a_2,\dots$ различных чисел (в нашем случае из $\mathbb Z_{\ge 0}$)
$$
h(t)=h(a_0)+\sum_{i=1}^\infty\Delta_h[a_0,\dots,a_i]\prod_{j=0}^{i-1}(t-a_j),
$$
где
$$
\Delta_h[a_0,\dots,a_i]:=\sum_{j=0}^i\frac{h(a_j)}{\prod_{p=0,p\ne j}^{i}(a_j-a_p)}.      
$$

  Согласно \cite{4}, p. 548, для $p\in\mathbb C$
$$
Q_r(n,m):=\sum_{k\ge 1}\frac{Q_r(n,k,m)}{k^{r(p+m)}}=\frac{\zeta(rp+r-m)}{\zeta(rp+r)}\sigma_{-p}(n,r),\eqno (9)
$$
где
$$
\sigma_p(n,r):=\sum_{d^r|n}d^{rp}.          
$$
Докажем аналогичную формулу для многочленов $M_r(t;n,k)$.

\begin{prop}\label{prop4} 
$$
M_r(t;n):=\sum_{k\ge 1}\frac{M_r(t;n,k)}{k^{rp}}=\frac{Li_{rp+r}(t)}{\zeta(rp+r)}\sigma_{-p}(n,k).\eqno (9)
$$
\end{prop}

  Доказательство. Используя правило умножения рядов Дирихле (коэффициенты произведения являются свёртками Дирихле коэффициентов сомножителей) и равенство \cite{4}, p. 548,
$$
\sum_{k\ge 1}\frac{c_r(n,k)}{k^{rp}}=\frac{\sigma_{1-p}(n,r)}{\zeta(rp)}
$$
получим
$$
\sum_{k\ge 1}\frac{M_r(t;n,k)}{k^{rp}}=\sum_{k\ge 1}\frac{1}{k^{rp+r}}\sum_{d|k}c_r(n,k/d)t^d
$$
$$
=\sum_{k\ge 1}\frac{c_r(n,k)}{k^{rp+r}}\sum_{k\ge 1}\frac{t^k}{k^{rp+r}}=\frac{Li_{rp+r}(t)}{\zeta(rp+r)}\sigma_{-p}(n,k).\blacksquare
$$

  Для производящих рядов Дирихле $Q_r(n,m)$ и $M_r(t;n)$ Предложение 3 примет следующий вид (полагаем $M(t;n):=M_1(t;n)$ и $Q(n,m):=Q_1(n,m)$).

\begin{prop}\label{prop3} 
$$
\left.\left(\frac{d}{dt}\right)^lM_r(t;n)\right|_{t=1}=\sum_{m=1}^ls(l,m)Q_r(n,m),\quad l\ge 1,
$$
$$
M_r(t;n)=\sum_{k:k^r|n}\frac{1}{k^{rp}}+\sum_{l\ge 1}\frac{(t-1)^l}{l!}\sum_{m=1}^ls(l,m)Q_r(n,m),
$$
$$
M_r(e^\lambda;n)=\sum_{k:k^r|n}\frac{1}{k^{rp}}+\sum_{l=1}^{\infty}\frac{\lambda^l}{l!}Q_r(n,l),
$$
в частности, для $r=1$
$$
\left.\left(\frac{d}{dt}\right)^lM(t;n)\right|_{t=1}=\sum_{m=1}^ls(l,m)Q(n,m),\quad l\ge 1,
$$
$$
M(t;n)=\sum_{k:k|n}\frac{1}{k^p}+\sum_{l\ge 1}\frac{(t-1)^l}{l!}\sum_{m=1}^ls(l,m)Q(n,m),
$$
$$
M(e^\lambda;n)=\sum_{k:k|n}\frac{1}{k^p}+\sum_{l=1}^{\infty}\frac{\lambda^l}{l!}Q(n,l).\blacksquare
$$
\end{prop}

\end {document}